\documentclass[10pt]{article}
\usepackage{latexsym}
\usepackage{amsthm}
\usepackage{amssymb}
\newtheorem{theorem}{Theorem}
\newtheorem{definition}{Definition}
\newtheorem{proposition}{Proposition}

\newtheorem{corollary}{Corollary}
\title{$k$-hypermonogenic automorphic forms}
\author{
D.~Constales \thanks{Department of
Mathematical Analysis, Ghent University, Building
S-22, Galglaan 2, B-9000 Ghent, Belgium. Financial support from BOF/GOA 01GA0405 of Ghent University gratefully acknowledged.  E-mail:
{\tt dc@cage.UGent.be}}
\and R. S. Krau{\ss}har \thanks{Department of
Mathematical Analysis, Ghent University, Building
S-22, Galglaan 2, B-9000 Ghent, Belgium. E-mail:
{\tt krauss@cage.UGent.be}}
\and John Ryan \thanks{Department of Mathematics, University of Arkansas, Fayetteville, AR 72701, USA {\tt jryan@uark.edu}}}
\begin{document}
\maketitle
\date
\begin{abstract}
In this paper we deal with monogenic and $k$-hypermonogenic automorphic forms on arithmetic subgroups of the Ahlfors-Vahlen group. Monogenic automorphic forms are exactly the $0$-hypermonogenic automorphic forms.\\ 
In the first part we establish an explicit relation between $k$-hypermonogenic automorphic forms and  Maa{\ss} wave forms. In particular, we show how one can construct from any arbitrary non-vanishing monogenic automorphic form a Clifford algebra valued Maa{\ss} wave form.\\ 
In the second part of the paper we compute the Fourier expansion of the $k$-hypermonogenic Eisenstein series which provide us with the simplest non-vanishing  examples of $k$-hypermonogenic automorphic forms.  
\end{abstract}
{\bf AMS-Classification}: 11F03, 11F36, 11F30, 11F37, 30G35
\section{Introduction}

The theory of higher dimensional Maa{\ss} wave forms has become a major topic of study in analytic number theory. Maa{\ss} wave forms are automorphic forms that are complex-valued eigensolutions to the Laplace-Beltrami operator 
\begin{equation}
\label{LB}
\Delta_{LB} = x_n^2 (\sum_{i=0}^n\frac{\partial^2 }{\partial x_n^2}) - (n-1) x_n \frac{\partial }{\partial x_n}.
\end{equation}
The classical setting is $n+1$-dimensional upper half-space in the framework of the action of discrete arithmetic subgroups of the orthogonal group. Its study was initiated in 1949 by H. Maa{\ss} in \cite{Maa1,Maa2}. In the late 1980s this study had a major boost by breakthrough works of J.~Elstrodt, F.~Grunewald and J.~Mennicke, \cite{EGM85,EGM88,EGM90}, A.~Krieg \cite{Kri2,Kri88}, V. Gritsenko \cite{Grit} among many others.  

\par\medskip\par 
 
Recently, in~\cite{KraHabil} another class of automorphic forms on these arithmetic groups has been considered. The context is again $n+1$-dimensional upper half-space. The classes of automorphic forms considered in \cite{KraHabil}, however,  have different analytic and mapping properties. They are null-solutions to iterates of the Euclidean Dirac  operator $D:=\sum_{i=1}^{n+1} \frac{\partial }{\partial x_i} e_i$  and, in general, they take values in real Clifford algebras.  These are called $k$-monogenic automorphic forms.  The monogenic automorphic forms on upper half-space in turn can be embedded into the general framework of $k$-hypermonogenic automorphic forms. The class of $k$-hypermonogenic functions \cite{Leutwiler, slberlin, slnewnew} contains the 1-monogenic functions as the particular  case $k=0$.  $1$-monogenic functions are more known simply as monogenic functions, \cite{DSS}.

\par\medskip\par

An important question that remained open in this context is to understand whether these classes of automorphic forms are linked with Maa{\ss} wave forms. 

\par\medskip\par 

In Section~3 we shed some light on this problem after having introduced the preliminary tools. We develop explicit relations between the classes of $k$-hyper- monogenic automorphic forms and eigensolutions to the Laplace-Beltrami operator. In particular, we explain how we can construct Clifford algebra valued Maa{\ss} wave forms from monogenic automorphic forms. 
   
\par\medskip\par

In Section~4 we compute the Fourier expansion of $k$-hypermonogenic automorphic forms.  
The appearance of Bessel-K functions in its Fourier expansion reflects once more the close relation between the class of  $k$-hypermonogenic automorphic forms and Maa{\ss} wave forms, see \cite{Maa2,EGM85,Kri2,Kri88}. In particular, we compute the Fourier coefficients of the $k$-hypermonogenic Eisenstein series explicitly.  

\par\medskip\par

The appearance of variants of the Epstein zeta function provides a further nice analogy to the classical theory of complex analytic automorphic forms. See for example \cite{Freitag,KK, schoeneberg,Terras}.  

\section{Preliminaries}
\subsection{Clifford algebras}
We introduce the basic notions of real Clifford algebras over the Euclidean space $\mathbb{R}^n$. For details, see for instance~\cite{DSS}. Throughout this paper, $\{e_1,\ldots,e_n\}$ stands for the standard orthonormal basis in the Euclidean space $\mathbb{R}^n$ and $Cl_n$ denotes its associated real Clifford algebra
in which the relation $e_i e_j + e_j e_i = - 2 \delta_{ij}$ holds.  A basis for
the algebra $Cl_{n}$ is given by $1, e_{1},\ldots,e_{n},\ldots,e_{j_{1}}\ldots e_{j_{r}},\ldots, e_{1}\ldots e_{n}$ 
where $1 \le j_{1}<\ldots j_{r} \le n$. Every element $a$ from $Cl_n$ can hence be written in the form 
$a = \sum_{A \subseteq P(\{1,\ldots,n\})} a_A e_A$ where $a_A$ are uniquely defined real numbers.  Here, $e_{\emptyset} = e_0 = 1$. 
Elements from $Cl_n$ of the reduced form $x_{0} + x_1 e_1 + \cdots + x_n e_n$ are called paravectors. $x_{0}$ is called the scalar part of $x$ and will be denoted by $Sc(x)$. In the case $n=1$ the Clifford algebra $Cl_1$ only consists of paravectors of the form $x_0 + x_1 e_1$ which in turn can be identified with the complex numbers. In the case $n=2$ we have $Cl_2 \stackrel{\sim}{=} \mathbb{H}$, $\mathbb{H}$ standing for the Hamiltonian quaternions. These in turn can be identified with the set of paravectors in $\mathbb{R} \oplus \mathbb{R}^3$ which have the form $x=x_0 + x_1 e_1 + \cdots + x_3 e_3$ when we identify the particular imaginary unit $e_3$ with the product $e_1 e_2$. 

The Clifford conjugation is defined by $\overline{ab} = \overline{b}\;\overline{a}$ where $\overline{e_j} = -e_j$ for all $j=1,2,\ldots,n$.  Each paravector $x \in \mathbb{R} \oplus \mathbb{R}^n \backslash\{0\}$ has an inverse element, given by $\overline{x}/\|x\|^2$. 
 
The reversion anti-automorphism is defined by  $\widetilde{ab} = \tilde{b} \tilde{a}$, where $\widetilde{e_{j_{1}}\ldots e_{j_{r}}} = e_{j_{r}}\ldots e_{j_{1}}$.
The main involution is defined by $(ab)^{'} = a^{'}b^{'}$ and $e_i^{'}=-1$ for $i=1,\ldots,n$.
Furthermore, we also need the following automorphism $*: Cl_n \rightarrow Cl_n$ defined by the
relations: $e_n^{*}=-e_n$, $e_i^{*}= e_i$ for $i=0,1,\ldots,n-1$ and $(ab)^{*}=a^{*}b^{*}$. Any element $a \in Cl_n$ may be uniquely decomposed in the form $a = b + c e_n$, where $b$ and $c$ belong to $Cl_{n-1}$. Based on this definition one defines the mappings $P: Cl_n \rightarrow Cl_{n-1}$ and $Q: Cl_n \rightarrow Cl_{n-1}$ by $Pa=b$ and $Qa=c$.

\subsection{Dirac type operators over $n+1$-dimensional Euclidean and hyperbolic spaces}
{\bf Monogenic functions.} Let $U \subseteq \mathbb{R} \oplus \mathbb{R}^n$ be an open set. Then a real differentiable function $f: U \rightarrow Cl_n$ that  satisfies $Df = 0$, resp. $fD = 0$, where $D := \frac{\partial }{\partial x_0}  + \frac{\partial }{\partial x_1} e_1 + \cdots + \frac{\partial }{\partial x_n}e_n$
is the Euclidean Cauchy-Riemann operator, is called left (resp. right) monogenic or left (resp. right) Clifford-holomorphic. Due to the non-commutativity of $Cl_{n}$ for $n>1$, both
classes of functions do not coincide with each other. However, $f$ is left monogenic, if and only if $\tilde{f}$ is right monogenic. 

The Cauchy-Riemann operator factorizes the Euclidean Laplacian $\Delta = \sum_{j=0}^n \frac{\partial^2}{\partial x_j^2}$, viz $D \overline{D} = \Delta$. Every real component of a monogenic function is hence harmonic.  

\par\medskip\par

The $D$-operator is left quasi-invariant under M\"obius transformations. Following for example \cite{a,EGM87}, M\"obius transformations in $\mathbb{R} \oplus \mathbb{R}^n$ can be
represented as
$$
T: \mathbb{R} \oplus \mathbb{R}^n \cup \{\infty\} \rightarrow \mathbb{R} \oplus \mathbb{R}^n \cup \{\infty\},\; T(x)=(ax+b)(cx+d)^{-1}
$$
with coefficients $a,b,c,d$ from $Cl_n$ that can all be written as products of paravectors from $\mathbb{R} \oplus\mathbb{R}^n$ and that satisfy $a\tilde{d}-b\tilde{c}\in \mathbb{R}\backslash\{0\}$ and $a^{-1}b, c^{-1}d \in \mathbb{R} \oplus\mathbb{R}^n$ if $c \neq 0$ or $a \neq 0$, respectively. These conditions are often called Ahlfors-Vahlen conditions. 

The set that consists of Clifford valued matrices $\left(\begin{array}{cc} a & b \\ c & d \end{array}\right)$ whose coefficients satisfy the above conditions forms a group under matrix multiplication and it is called the general Ahlfors-Vahlen group, $GAV(\mathbb{R} \oplus\mathbb{R}^n)$. The subgroup consisting of those matrices from $GAV(\mathbb{R} \oplus\mathbb{R}^n)$ that  satisfy $a\tilde{d}-b\tilde{c}=1$ is called the special Ahlfors-Vahlen group. It is denoted by $SAV(\mathbb{R} \oplus\mathbb{R}^n)$. 

\par\medskip\par

Now assume that $M \in GAV(\mathbb{R} \oplus\mathbb{R}^n)$. If $f$ is a left monogenic function in the  variable
$y=M<x>=(ax+b)(cx+d)^{-1}$, then the function $\frac{\overline{cx+d}}{\|cx+d\|^{n+1}}f(M<x>)$ is again left monogenic in the variable $x$, cf. \cite{r85}.

\par\medskip\par

{\bf $k$-hypermonogenic functions.}
The class of monogenic functions belongs to the more general class of so-called $k$-hypermonogenic functions, introduced in \cite{slberlin}. These are defined as the null-solutions to the system
\begin{equation}\label{kh}
D f + k \frac{(Q f)^{'}}{x_n} = 0
\end{equation}
where $k \in \mathbb{R}$. 

\par\medskip\par

Under a M\"{o}bius transformation $y=M<x>=(ax+b)(cx+d)^{-1}$ such a solution $f(y)$ is transformed to the $k$-hypermonogenic function
\begin{equation}
\label{kinv}
F(x):=\frac{\overline{cx+d}}{\|cx+d\|^{n+1-k}}f(M<x>).
\end{equation}
See for instance \cite{slnewnew}.
In the case $k=0$, we are dealing with the set of left monogenic functions introduced earlier. The particular solutions associated to the case $k=n-1$ coincide with the null-solutions to the hyperbolic Hodge-Dirac operator with respect to the hyperbolic metric on upper half space. These are often called hyperbolic monogenic functions or simply hypermonogenic functions, see \cite{Leutwiler}.

\subsection{Discrete arithmetic subgroups of $GAV(\mathbb{R} \oplus\mathbb{R}^{n})$}
In this paper we deal with null-solutions to the above systems that are quasi-invariant under arithmetic subgroups of the general Ahlfors-Vahlen group that act totally discontinuously on upper half-space 
$$H^{+}(\mathbb{R} \oplus\mathbb{R}^n) =\{x = x_0 + x_1 e_1 + \cdots + x_n e_n \in \mathbb{R} \oplus\mathbb{R}^n\;:\;x_n > 0\}.$$ 
These can be regarded as generalizations of the classical holomorphic modular forms in the context of monogenic and $k$-hypermonogenic functions.  

\par\medskip\par
The automorphism group of upper half-space $H^{+}(\mathbb{R} \oplus\mathbb{R}^n)$ is the group $SAV(\mathbb{R} \oplus\mathbb{R}^{n-1})$.   
Arithmetic subgroups of the special Ahlfors-Vahlen group that act totally discontinuously on upper half-space are for instance considered in \cite{Maa2,EGM87,EGM90}. For convenience, we first recall
the definition of the rational Ahlfors-Vahlen group on $H^+(\mathbb{R} \oplus\mathbb{R}^n)$.
\begin{definition}
The rational Ahlfors-Vahlen group $SAV(\mathbb{R} \oplus\mathbb{R}^{n-1},\mathbb{Q})$ is the set of matrices
$\left(\begin{array}{cc} a & b \\ c & d \end{array}\right)$ from $SAV(\mathbb{R} \oplus\mathbb{R}^{n-1})$ that satisfy \\[0.2cm]
$ (i)\quad a \overline{a}, b \overline{b}, c \overline{c}, d \overline{d} \in {\mathbb{Q}},\\
  (ii) \quad a\overline{c},b\overline{d} \in {\mathbb{Q}}
\oplus \mathbb{Q}^{n},\\
  (iii)  \quad a x \overline{b} + b \overline{x} \;\overline{a}, c x \overline{d} + d \overline{x} \;\overline{c} \in    {\mathbb{Q}}\quad (\forall x \in \mathbb{Q} \oplus 
\mathbb{Q}^{n}),\\
 (iv) \quad a x \overline{d}+b \overline{x}\;\overline{c}  
\in {\mathbb{Q}}
\oplus \mathbb{Q}^{n}
\quad (\forall x \in \mathbb{Q} \oplus \mathbb{Q}^{n})$.
\end{definition}
Next we need
\begin{definition}
A $\mathbb{Z}$-order in a rational Clifford algebra is a subring $R$ such that the additive group of $R$ is finitely generated and contains a $\mathbb{Q}$-basis of the
Clifford algebra.
\end{definition}
The following definition, cf. \cite{EGM90}, provides us with a whole class of arithmetic subgroups  of the Ahlfors-Vahlen group which act totally discontinuously on upper half-space. 
\begin{definition}
Suppose that $p \le n-1$. 
Let ${\cal{I}}$ be a ${\mathbb{Z}}$-order in ${Cl}_{n}$ which is stable under the reversion and the main involution ${'}$ of ${Cl}_{n}$. Then
$$
\Gamma_{p}({\cal{I}}) := SAV(\mathbb{R} \oplus\mathbb{R}^{p},{\mathbb{Q}}) \cap Mat(2,{\cal{I}}).
$$
For an $N \in {\mathbb{N}}$ the principal congruence subgroup of
$\Gamma_{p}({\cal{I}})$ of level $N$ is defined by
$$
\Gamma_{p}({\cal{I}})[N] := \Big\{ \left(\begin{array}{cc} a & b\\ c & d \end{array}\right) \in \Gamma_{p}({\cal{I}}) \;\Big|\; a-1,b,c,d-1 \in N {\cal{I}}\Big\}.
$$
\end{definition}
Notice that all the groups $\Gamma_{p}({\cal{I}})[N]$ have finite index in $\Gamma_{p}({\cal{I}})$. Therefore, all of them are discrete groups and act totally discontinuously on upper half-space $H^{+}(\mathbb{R} \oplus{\mathbb{R}}^n)$ provided $p \le n-1$. 
The proof of the total discontinuous action can be done in the same way as in \cite{Freitag}. 

\par\medskip\par
{\bf Special notation}. 
In the sequel let us denote the subgroup of translation matrices contained in $\Gamma_{p}({\cal{I}})[N]$ by ${\cal{T}}_{p}({\cal{I}})[N]$. The associated $p+1$-dimensional lattice in $\mathbb{R} \oplus \mathbb{R}^{p}$ will be denoted by $\Lambda_{p}({\cal{I}})[N]$ and its dual by $\Lambda_{p}^{*}({\cal{I}})[N]$. The dual lattice is contained in the subspace $\mathbb{R} \oplus \mathbb{R}^{p}$, too.  

\par\medskip\par

The simplest concrete examples for $\Gamma_{p}({\cal{I}})$ are obtained by taking for ${\cal{I}}$ the standard ${\mathbb{Z}}$-order in the Clifford algebras $Cl_p$, i.e.,  ${\cal{O}}_p := \sum_{A \subseteq  P(1,\ldots,p)} \mathbb{Z} e_A$ where $p \le  n-1$.  
In this case, the group $\Gamma_{p}({\cal{I}})$ coincides with the special hypercomplex modular group considered in \cite{KraHabil}. This one is generated by the matrices 
$$
J:= \left(\begin{array}{cc} 0 & -1 \\ 1 & 0 \end{array}\right),
T_{1} := \left(\begin{array}{cc} 1 & 1 \\ 0 & 1  \end{array}\right), 
T_{e_1} := \left(\begin{array}{cc} 1 & e_1 \\ 0 & 1  \end{array}\right),\ldots, T_{e_p}:= \left(\begin{array}{cc} 1 & e_p \\ 0 & 1  \end{array}\right).
$$
In this particular case, we have ${\cal{T}}_{p}({\cal{O}}_p) = <T_1, T_{e_1},\ldots,T_{e_p}>$ and the associated period lattice is the orthonormal $p+1$-dimensional lattice $\mathbb{Z} + \mathbb{Z} e_1 + \cdots + \mathbb{Z} e_p$ which is self-dual. Its standard fundamental period cell is $[0,1]^{p+1}$. 

\par\medskip\par
In the case $n=3$ the set $H^{+}(\mathbb{R} \oplus \mathbb{R}^3)$  can be identified with upper quaternionic half-space. In this setting further important examples of ${\cal{I}}$ are the quaternionic orders, in particular the Hurwitz order $\mathbb{Z} + \mathbb{Z} e_1  +\mathbb{Z} e_2 + \frac{1}{2} \mathbb{Z} (1+e_1 + e_2+e_3)$,  considered for example in \cite{H,Kri2,Kri88,Reid}. In the quaternionic context the third unit $e_3$ is identified with $e_1 e_2$.

\subsection{$k$-hypermonogenic automorphic forms}

As already mentioned, any M\"obius transformation $T(x)=M<x>$ induced by  matrices $M$ belonging to $SAV(\mathbb{R} \oplus \mathbb{R}^{n-1})$ (as well as any of its subgroups) preserves upper half-space. Therefore, in view of~(\ref{kinv}), if $f$ is a
function that is left $k$-hypermonogenic on the whole half-space then so is $\frac{\overline{cx+d}}{\|cx+d\|^{n+1-k}} f(M<x>)$. 

Next we introduce 
\begin{definition}
Let $p \le n-1$. 
A left $k$-hypermonogenic function $f: H^{+}(\mathbb{R} \oplus\mathbb{R}^n) \rightarrow Cl_n$ is called a left $k$-hypermonogenic automorphic form on $\Gamma_{p}({\cal{I}})[N]$ if for all $x \in H^{+}(\mathbb{R} \oplus \mathbb{R}^n)$  
\begin{equation}
\label{stern}
f(x)=\frac{\overline{cx+d}}{\|cx+d\|^{n+1-k}} f(M<x>) 
\end{equation}
for all $M \in \Gamma_{p}({\cal{I}})[N]$. 
\end{definition}
In the case $k=0$ we re-obtain the class of left monogenic automorphic forms described in \cite{KraHabil}. 
As the following proposition shows, there is a direct relation between $k$-hypermonogenic automorphic forms and  $-k$-hypermonogenic automorphic forms.  
\begin{proposition}\label{-k}
Suppose that $k\in \mathbb{R}$ and that $p$ is a positive integer with $p < n$. If $f:H^{+}(\mathbb{R} \oplus\mathbb{R}^n) \rightarrow Cl_n$ is a $k$-hypermonogenic automorphic form on $\Gamma_{p}({\cal{I}})[N]$ satisfying (\ref{stern}), then the function $g:H^{+}(\mathbb{R} \oplus \mathbb{R}^n) \rightarrow Cl_n$ defined by $g(x) := \frac{f(x)e_n}{x_n^k}$ is $-k$ hypermonogenic and satisfies 
\begin{equation}
g(x) =  \frac{\overline{cx+d}}{\|cx+d\|^{n+1+k}} g(M<x>)
\end{equation}
for all $M \in \Gamma_{p}({\cal{I}})[N]$. 
\end{proposition}
\begin{proof}
Following for example~\cite{slnew}, a function $f$ is $k$-hypermonogenic, if and only if $\frac{f e_n}{x_n^k}$ is $-k$  hypermonogenic.  
Suppose that $f$ is a $k$-hypermonogenic automorphic form satisfying (\ref{stern}). Then  
$g(x) := \frac{f(x)e_n}{x_n^k}$ 
is $-k$-hypermonogenic and satisfies the transformation law 
\begin{eqnarray*}
g(M<x>) &=& \frac{1}{(M<x>)_n^k} f(M<x>) e_n\\
&=& \frac{\|cx+d\|^{2k}}{x_n^k} (\|cx+d\|^{n+1-k}) (\overline{cx+d})^{-1} f(x) e_n\\
&=& \|cx+d\|^{n+1+k}(\overline{cx+d})^{-1} g(x).
\end{eqnarray*}
\end{proof}
Therefore, it is sufficient to restrict to non-positive values of $k$ in all that follows.

\par\medskip\par
  
For each ${\cal{I}}$ there exists a minimal positive integer $N_0({\cal{I}})$ such that neither the negative identity matrix $-I$ nor the other diagonal matrices of the form $\left(\begin{array}{cc} \tilde e_A& 0 \\ 0 & e_A^{-1} \end{array} \right)$ where $A \subseteq P(1,\ldots,p)$ are no longer included in all principal congruence subgroups $\Gamma_{p}({\cal{I}})[N]$ with $N \ge N_0({\cal{I}})$. In the case where ${\cal{I}}$ is the standard $\mathbb{Z}$-order ${\cal{O}}_p$, we have $N_0({\cal{O}}_p) = 3$, see \cite{KraHabil}. For all $N < N_0({\cal{I}})$, only the zero function satisfies (\ref{stern}). However, for all $N \ge N_0({\cal{I}})$, one can construct non-trivial $k$-hypermonogenic automorphic forms that have the transformation behavior (\ref{stern}). 

\par\medskip\par

The simplest non-trivial examples of $k$-hypermonogenic automorphic forms on the groups $\Gamma_p({\cal{I}})[N]$ with $N \ge N_0({\cal{I}})$ are the $k$-hypermonogenic generalized Eisenstein series:

\begin{definition}
Let $N \ge N_0({\cal{I}})$. For $p < n$ and real $k$ with $k < n-p-2$ the $k$-hypermonogenic Eisenstein series on the group $\Gamma_{n-1}({\cal{I}})[N]$ acting on $H^{+}(\mathbb{R} \oplus\mathbb{R}^n)$ are defined by
\begin{equation}\label{khypereisen}
\varepsilon_{k,p,N}(x) := \sum\limits_{M:{\cal{T}}_{p}({\cal{I}})[N] \backslash \Gamma_{p}({\cal{I}})[N]} \frac{\overline{cx+d}}{\|cx+d\|^{n+1-k}}.
\end{equation}
\end{definition}
These series converge for $k < n-p-2$ absolutely and uniformly on each compact subset of $H^{+}(\mathbb{R} \oplus\mathbb{R}^n)$. A majorant is 
$$
\sum\limits_{M: {\cal{T}}_p({\cal{I}})[N] \backslash \Gamma_p({\cal{I}})[N]} \frac{1}{\|c e_n +d \|^{\alpha}}
$$
whose absolute convergent abscissa is $\alpha > p+2$, cf. for example \cite{KraHabil,EGM90}. For $p=n-1$, this majorant converges absolutely for all $k <-1$.  

\par\medskip\par

The non-vanishing behavior for $N \ge N_0({\cal{I}})$ can easily be established by considering the limit $\lim_{x_n \rightarrow +\infty} \varepsilon_{k,p,M}(x_n e_n)$ which equals $+1$ in these cases,~cf.~\cite{KraBMS}.  

\par\medskip\par 

{\bf Remarks.}  
In the case $k=0$ and ${\cal{I}} = {\cal{O}}_p$ the series (\ref{khypereisen}) coincide with the monogenic Eisenstein series considered in \cite{KraHabil,KraBMS} in the cases $p < n-2$.  

\par\medskip\par

By adapting the Hecke trick from \cite{Freitag} can also introduce Eisenstein series of lower weight. This is shown in Section 4 of \cite{paper4}. In particular, for $N \ge N_0({\cal{I}})$ the series 
\begin{equation}\label{0hypereisen}
\varepsilon_{0,N-1,M}(x) := \lim\limits_{s \rightarrow 0^+}\sum\limits_{M:{\cal{T}}_{n-1}({\cal{I}})[N] \backslash \Gamma_{n-1}({\cal{I}})[N]} \Bigg(\frac{x_n}{\|cx+d\|^2}\Bigg)^s \frac{\overline{cx+d}}{\|cx+d\|^{n+1}}.
\end{equation}
defines a well-defined non-vanishing left monogenic Eisenstein series on the groups $\Gamma_{n-1}({\cal{I}})[N]$ in  upper half-space variable $x$, cf.~\cite{paper4}~Section~4. In \cite{paper4}~Section~4 this is done for the particular case ${\cal{I}} = {\cal{O}}_p$. However, the transition to the context of more general orders ${\cal{I}}$ follows identically along the same lines.   

\par\medskip\par

In view of Proposition~\ref{-k}, we can directly construct non-vanishing $j$-hypermonogenic Eisenstein series for positive $j$ from the $k$-hypermonogenic Eisenstein series of negative $k$, simply by forming
$$
E_{-k,p,N} (x) := \frac{\varepsilon_{k,p,N}(x) e_n}{x_n^k}
$$
which then satisfy the transformation law
$$
E_{-k,p,N} (x) := \frac{\overline{cx+d}}{\|cx+d\|^{n+1+k}} E_{-k,p,N}(M<x>)
$$    
for all $M \in \Gamma_p({\cal{I}})[N]$.  

\section{Relation to Maa{\ss} wave forms} 

As one directly sees, the theory of monogenic automorphic forms fits as a special case within the general framework of  $k$-hypermonogenic automorphic forms. We can say more:

\par\medskip\par

Suppose that $f = Pf + Qf e_n$ is a $k$-hypermonogenic function, where $k$ is some arbitrary fixed real number. Then, 
following for example \cite{slberlin} all the real components of $Pf$ are contained in the set of $k$-hyperbolic harmonic functions. These are solutions to the system
\begin{equation}
\label{khypharm}
x_n \Delta u- k \frac{\partial u}{\partial x_n} = 0.
\end{equation}
These solutions are also quasi-invariant under M\"obius transformations that act on upper half-space: If $f$ is a solution to (\ref{khypharm}), then following \cite{slnew},
\begin{equation}
\label{harminv}
F(x) = \frac{1}{\|cx+d\|^{n-k-1}} f(M<x>)
\end{equation}
is $k$-hyperbolic harmonic, too. 
Notice that in the particular case $k=n-1$ the correction factor disappears. This attributes a special role to the function  class of $(n-1)$-hypermonogenic functions.  

\par\medskip\par

The solutions to (\ref{khypharm}) in turn are directly related to the Maa{\ss} wave equation. Following for example \cite{Leutw87}, if $u$ is a solution to (\ref{khypharm}) then  $g(x) = x_n^{-(1-n+k)/2}u(x)$ is a solution to 
\begin{equation}
\label{maass}
\Delta g - \frac{n-1}{x_n} \frac{\partial g}{\partial x_n} + \lambda \frac{g}{x_n^2} = 0
\end{equation}
where $\lambda = \frac{1}{4}(n^2 -(k+1)^2)$. 
The solutions to (\ref{maass}) have the property that they are directly preserved by all M\"obius transformations that  act on upper half-space. Each solution $g$ is an eigensolution to the Laplace-Beltrami operator (\ref{LB}) associated to the  fixed eigenvalue $-\frac{1}{4}(n^2-(k+1)^2).$ 

\par\medskip\par

Let $p < n$. Now suppose that $f: H^{+}(\mathbb{R} \oplus\mathbb{R}^n) \rightarrow Cl_n$ is a $k$-hypermonogenic automorphic form on $\Gamma_{p}({\cal{I}})[N]$ of weight $(n-k)$, satisfying the transformation law $f(x) = \frac{\overline{cx+d}}{\|cx+d\|^{n+1-k}} f(M<x>)$ for all $M \in \Gamma_{p}({\cal{I}})[N]$.
 
Since $(M<x>)_n = \frac{x_n}{\|cx+d\|^2}$, the function $g(x) = x_n^{-(1-n+k)/2}f(x)$ thus satisfies for all $M \in \Gamma_{p}({\cal{I}})[N]$:
\begin{eqnarray*}
g(M<x>) &=& {(M<x>)_n}^{-\frac{1-n+k}{2}} f(M<x>) \\
&=& \Bigg(\frac{x_n}{\|cx+d\|^2}\Bigg)^{-\frac{1-n+k}{2}} \|cx+d\|^{n+1-k} (\overline{cx+d})^{-1} f(x)\\
&=& x_n^{-\frac{1-n+k}{2}} \frac{1}{\|cx+d\|^{n-k-1}} \|cx+d\|^{n+1-k} (\overline{cx+d})^{-1} x_n^{\frac{1-n+k}{2}} g(x). 
\end{eqnarray*}
Hence, 
$$
g(x) = \frac{\overline{cx+d}}{\|cx+d\|^2} g(M<x>). 
$$
Unfortunately, if $f = Pf + Qf e_n$ is a $k$-hypermonogenic automorphic form with respect to $\Gamma_{p}({\cal{I}})[N]$, then $Pf$ is in general not an automorphic form with respect to the full group $\Gamma_p[N]({\cal{I}})$. However, only the components of the $P$-part of $f$ satisfy (\ref{khypharm}). The associated function $g(x)= x_n^{-(1-n+k)/2}f(x)$ is exactly a $\Gamma_p({\cal{I}})[N]$-invariant eigenfunction to the Laplace-Beltrami operator for the eigenvalue  $-\frac{1}{4}(n^2-(k+1)^2)$, if $Q f= 0$. In this case we are dealing with monogenic functions. Actually, we can construct from every left monogenic automorphic form a Maa{\ss} wave form. Following \cite{slberlin}, the real components of the $Q$-part of a $k$-hypermonogenic function satisfies the equation
$$
x_n^2 \Delta u - k x_n \frac{\partial u}{\partial x_n} + k u = 0.
$$
In particular, if $f$ is monogenic, then all the real components of the $P$-part and of the  $Q$-part satisfy  
$$
x_n^2 \Delta u = 0.   
$$
Consequently, all real components of $g(x) := x_n^{-(1-n+k)/2}f(x)$ are eigenfunctions to the Laplace-Beltrami operator for the particular eigenvalue  $-\frac{1}{4}(n^2-1)$. 

We thus arrive at the following result:

\begin{theorem}\label{relation}
Suppose that $f$ is a left monogenic automorphic form on $\Gamma_{p}({\cal{I}})[N]$ of weight $n$, satisfying 
$f(x) = \frac{\overline{cx+d}}{\|cx+d\|^{n+1}} f(M<x>)$ for all $M \in \Gamma_{p}({\cal{I}})[N]$.  Then $g(x) = x_n^{-(1-n)/2}f(x)$ is a quasi-$\Gamma_p({\cal{I}})[N]$-invariant Maa{\ss} wave form associated to the fixed eigenvalue $-\frac{1}{4}(n^2-1)$ and has the $-1$-weight automorphy factor $\frac{\overline{cx+d}}{\|cx+d\|^2}$.   
\end{theorem}
Notice that in the complex case $n=1$, the function $g$ is a null-solution to the Laplace-Beltrami operator. 

\par\medskip\par

Theorem~\ref{relation} reflects a connection between the class of monogenic automorphic forms  and the particular family of non-analytic automorphic forms on the Ahlfors-Vahlen group considered for instance by A. Krieg, J. Elstrodt et al. and others (see for example \cite{Kri88,EGM90,Kri03}). The non-analytic automorphic forms considered in \cite{Kri88,EGM90,Kri03}) are scalar-valued eigenfunctions of the Laplace-Beltrami operator 
associated to a special continuous spectrum of eigenvalues described in \cite{EGM90}. They are all totally invariant under the group action. 

\par\medskip\par

Notice that the Eisenstein series considered in this paper are Clifford algebra valued in general. They are associated to a fixed eigenvalue.

\par\medskip\par

Let us now analyze more into detail how the particular examples of non-analytic Maa{\ss} wave forms discussed in the above mentioned papers fit within the framework of $k$-hypermonogenic automorphic forms. 
First we observe, that evidently in all the other cases where $f$ is a $k$-hypermonogenic automorphic form on $\Gamma_{p}({\cal{I}})[N]$ of weight $(n-k)$ with $k \neq 0$ and $Qf \neq 0$, the $P$-part of the associated function $g(x) = x_n^{-(1-n)/2}f(x)$ is an eigenfunction of the Laplace-Beltrami operator for the eigenvalue $-\frac{1}{4}(n^2-(k+1)^2)$  
and satisfies 
$$
P\Big[g(x)\Big] = P \Big[\frac{\overline{cx+d}}{\|cx+d\|^2} g(M<x>)\Big]\;\; \forall M \in \Gamma_p({\cal{I}})[N].   
$$
However, a relation of the form $$P[j(M,x) g(M<x>)] = j(M,x) P[g(M<x>)]$$ only holds if $j(M,x)$ is a scalar-valued automorphy factor.    
On the other hand, as already mentioned, the real components of the $P$-part of a $k$-hypermonogenic automorphic are $k$-hyperbolic harmonic. 

Let us now suppose that $F$ is a $k$-hyperbolic harmonic automorphic form that has the following invariance behavior $$F(x) = \frac{1}{\|cx+d\|^{(n-k-1)}} F(M<x>),$$
dealing with a scalar-valued automorphy factor. The associated function  $G(x)= x_n^{-(1-n+k)/2}F(x)$ then has the property of  being totally invariant under the group action of $\Gamma_{p}({\cal{I}})[N]$, i.e. $G(x) = G(M<x>)$ for all $M \in \Gamma_{p}({\cal{I}})[N]$.  
\par\medskip\par
Non-trivial examples of $k$-hyperbolic harmonic functions $F$ with this transformation behavior are the $k$-hyperharmonic Eisenstein series:
\begin{equation}
{\cal{E}}_{k,p,N}(x) := \sum\limits_{M:{\cal{T}}_{p}({\cal{I}})[N] \backslash \Gamma_{p}({\cal{I}})[N]} \frac{1}{\|cx+d\|^{n-k-1}}\;\;\; k < n-p-1.
\end{equation}
Indeed, when applying the transformation $G(x)=x_n^{-(1-n+k)/2}{\cal{E}}_{k,p,N}(x)$, then we obtain the totally invariant Eisenstein series 
\begin{equation}
{\cal{Q}}_{k,p,N}(x):=\sum\limits_{M:{\cal{T}}_{p}({\cal{I}})[N] \backslash \Gamma_{p}({\cal{I}})[N]} 
\Bigg(\frac{x_n}{\|cx+d\|^2}\Bigg)^{(n-k-1)/2}\;\;\; k < n-p-1.
\end{equation}   
In the case $p=n-1$, we then recover exactly the particular non-analytic Eisenstein series from \cite{EGM90} which is associated to the fixed eigenvalue $-\frac{1}{4}(n^2-(k+1)^2)$. 

In the particular complex case $n=1$ and ${\cal{I}}=\mathbb{Z}$, this series coincides with the complex non-analytic Eisenstein series from \cite{Maa1,Terras} that is associated to the particular eigenvalue $-\frac{1}{4}((1-(k+1)^2)$. In the subcase $n=2$ and ${\cal{I}} = \mathbb{Z}[e_1]$ we obtain the non-analytic Eisenstein series on Picard's group introduced in \cite{EGM85} that is associated to the fixed eigenvalue $-\frac{1}{4}((4-(k+1)^2)$.  

In the quaternionic setting where $n=3$ and where ${\cal{I}}$ is the Hurwitz order, we recover one of the non-analytic Eisenstein series from \cite{Kri88}.

\section{Fourier transforms}

\subsection{The Fourier expansion of $k$-hypermonogenic automorphic forms} 

In the first part of this section we determine the general structure of the Fourier expansion of $k$-hypermonogenic automorphic forms. Afterwards we then compute the Fourier coefficients of the $k$-hypermonogenic Eisenstein series that we introduced in (\ref{khypereisen}) explicitly.   
For technical reasons it is convenient to compute first the Fourier transform of a function that satisfies the reduced $k$-hypermonogenic equation 
\begin{equation}
\label{reduced}
{\cal{D}} f + \frac{k}{x_n} (Q f)^{'} = 0
\end{equation}
in reduced upper half-space $H^{+}(\mathbb{R}^n) := \{x_1 e_1 + \cdots + x_n e_n \in \mathbb{R}^n, x_n > 0\}$, where ${\cal{D}} := \sum\limits_{i=1}^n \frac{\partial }{\partial x_i} e_i$ is the Euclidean Dirac operator.  
We prove 
\begin{theorem}\label{reducedkf}
Suppose that $f:H^{+}(\mathbb{R}^n) \rightarrow Cl_n$ is a solution to (\ref{reduced}). 
Then its image under the Fourier transform in the directions $x_1,\ldots,x_{n-1}$ has the form 
$$
\alpha(\underline{\omega},x_n) = x_n^{\frac{k+1}{2}} \Bigg[K_{\frac{k+1}{2}}(\|\underline{\omega}\| x_n) - i e_n  \frac{\underline{\omega}}{\| \underline{\omega}\|} K_{\frac{k-1}{2}}(\|\underline{\omega}\| x_n)\Bigg] \alpha(\underline{\omega}) 
$$
if $\underline{\omega} \in \mathbb{R}^{n-1}\backslash \{0\}$ 
and, for $\underline{\omega} = \underline{0}$, 
$$
\alpha(\underline{\omega},x_n) a(\underline{0}) + \alpha(\underline{0}) x_n^k
$$
where $a(\underline{0})$ and all $\alpha(\underline{\omega})$ are well-defined Clifford numbers from the Clifford sub algebra $Cl_{n-1}$. 
\end{theorem}
\begin{proof} Applying the Fourier transform (in the first $(n-1)$ components) on the reduced $k$-hypermonogenic equation
$$
{\cal D} f + \frac{k}{x_n} (Q f)^{'} = 0
$$
leads to the following differential equation on the Fourier images 
$$
\alpha = (\alpha(\underline{\omega},x_n),  (\underline{\omega} \in \mathbb{R}^{n-1}, x_n > 0)
$$ 
of its solutions 
\begin{equation}\label{f1}
\Big(i \underline{\omega} \alpha + e_n \frac{\partial \alpha}{\partial x_n}\Big) + \frac{k}{x_n} (Q \alpha)^{'} = 0.
\end{equation}
Next we split $\alpha$ into its $P$-part and $Q$-part:
$$
\alpha = P \alpha + e_n (Q \alpha)^{'}.
$$
Applying this decomposition to each single term of (\ref{f1}) leads to 
\begin{eqnarray}
i \underline{\omega} \alpha & = & i \underline{\omega} (P \alpha) - e_n(i \underline{\omega} Q \alpha )^{'}\\
e_n \frac{\partial \alpha}{\partial x_n} &=& - \frac{\partial }{\partial x_n}(Q \alpha)^{'} + e_n \Big(  \frac{\partial }{\partial x_n}(P \alpha)^{'}  \Big)^{'}\\
\frac{k}{x_n} (Q \alpha)^{'} &=& (\frac{k}{x_n} Q \alpha)^{'}. 
\end{eqnarray}
Since $\alpha$ is a solution to (\ref{f1}), we thus obtain the following system of two equations when collecting all terms that include or that do not include the factor $e_n$, respectively:
\begin{eqnarray}
i \underline{\omega} (P \alpha) - \frac{\partial }{\partial x_n} (Q \alpha)^{'}+ \frac{k}{x_n}(Q \alpha)^{'} &=& 0\\
 i \underline{\omega} (Q \alpha) + \frac{\partial }{\partial x_n} (P \alpha)^{'} &=& 0. 
\end{eqnarray} 
This system is equivalent to 
\begin{eqnarray}
\label{f1'}
\frac{\partial }{\partial x_n} (Q \alpha)^{'} &=& i \underline{\omega} (P \alpha) + \frac{k}{x_n}(Q \alpha)^{'}\\
\frac{\partial }{\partial x_n} (P \alpha) &=&  i \underline{\omega} (Q \alpha)^{'}. \label{f2'}
\end{eqnarray}
To proceed, we next decompose $\alpha$ as follows
\begin{equation}
\label{dec}
\alpha = A \alpha + \underline{\omega} B \alpha
\end{equation} 
where $A$ and $B$ are appropriately chosen Clifford operators, projecting $\alpha$ down to that part of the Clifford algebra which intersects the subspace generated by $\underline{\omega}$ only at $0$. 

For the $P$-part and the $Q$-part we have 
\begin{eqnarray}
P \alpha &=& A P \alpha + \underline{\omega} B P \alpha \\
(Q \alpha)^{'} &=& A(Q \alpha)^{'} + \underline{\omega}B (Q \alpha)^{'}.
\end{eqnarray}
If we apply this decomposition to (\ref{f1'}), then we obtain 
\begin{eqnarray*}
\frac{\partial }{\partial x_n} A (Q \alpha)^{'} + \underline{\omega} \frac{\partial }{\partial x_n}B(Q \alpha)^{'} &=& i \underline{\omega}(A (P \alpha) + \underline{\omega} B(P \alpha)) \\
&+& \frac{k}{x_n} (A (Q\alpha)^{'} + \underline{\omega} B(Q \alpha)^{'}).
\end{eqnarray*} From this equation, we in turn recover the following system 
\begin{eqnarray}
\label{f3}
\frac{\partial }{\partial x_n} A (Q \alpha)^{'} &=& -i \|\underline{\omega}\|^2 B(P \alpha) + \frac{k}{x_n} A(Q \alpha)^{'}\\
\label{f4}
\frac{\partial }{\partial x_n} B(Q \alpha)^{'} &=& i A(P \alpha) + \frac{k}{x_n} B(Q \alpha)^{'}.
\end{eqnarray}
\par\medskip\par
Similarly, when applying the decomposition (\ref{dec}) to (\ref{f2'}), we obtain the equation 
\begin{eqnarray*}
\frac{\partial }{\partial x_n} A(P \alpha) + \underline{\omega} \frac{\partial }{\partial x_n} B(P \alpha) &=&  i \underline{\omega} (A(Q \alpha)^{'} + \underline{\omega} B(Q \alpha)^{'}) \\
&=&  i \underline{\omega} A (Q \alpha)^{'} - i\|\omega\|^2 B(Q \alpha)^{'}.
\end{eqnarray*}
This finally leads to the following system 
\begin{eqnarray}\label{f5}
\frac{\partial }{\partial x_n} A(P \alpha) &=& - i \|\omega\|^2 B (Q \alpha)^{'}\\
\frac{\partial }{\partial x_n} B(P \alpha) &=& i A (Q \alpha)^{'}. \label{f6}
\end{eqnarray}
For simplicity let us next use the abbreviations: $u_1 := A Q(\alpha)^{'}, v_1:= B P(\alpha)$ and $u_2:= AP(\alpha), v_2:= B Q(\alpha)^{'}$.  The set of equations (\ref{f3}) and (\ref{f6}) is a system involving only the functions $u_1$ and $v_1$:
\begin{eqnarray}
\frac{\partial }{\partial x_n} u_1 &=& - i \|\underline{\omega}\|^2 v_1 + \frac{k}{x_n} u_1\\
\frac{\partial }{\partial x_n} v_1 &=& i u_1.
\end{eqnarray}
The other set of equations (\ref{f4}) and (\ref{f5}) is a system involving only the other functions $u_2$ and $v_2$:
\begin{eqnarray}
\frac{\partial }{\partial x_n} u_2 &=& - i \|\omega\|^2 v_2 \\
\frac{\partial }{\partial x_n} v_2 &=& i u_2 + \frac{k}{x_n} v_2.
\end{eqnarray}
In the case where $\underline{\omega} \neq 0$,  
the general solutions to these systems turn out to be 
\begin{eqnarray}
u_1 &=& i x_n^{\frac{k+1}{2}} K_{\frac{k-1}{2}}(\|\underline{\omega}\| x_n) \alpha_1(\underline{\omega})\\
v_1 &=& \frac{1}{\|\underline{\omega}\|} x_n^{\frac{k+1}{2}} K_{\frac{k+1}{2}}(|\underline{\omega}| x_n)\alpha_1(\underline{\omega})\\
u_2 &=& x_n^{\frac{k+1}{2}} K_{\frac{k+1}{2}}(\|\underline{\omega}\| x_n)\alpha_2(\underline{\omega})\\
v_2 &=& -x_n^{\frac{k+1}{2}} \frac{i}{\|\underline{\omega}\|} K_{\frac{k-1}{2}}(\|\underline{\omega}\|x_n) \alpha_2(\underline{\omega})
\end{eqnarray}
where $\alpha_1(\underline{\omega})$ and $\alpha_2(\underline{\omega})$ may be Clifford constants that belong to the orthogonal complement of $span_{\mathbb{R}}(\underline{\omega},e_n)$. 

After re-substituting 
$$
\alpha = P(\alpha) + e_n Q(\alpha)^{'} = u_2 + \underline{\omega} v_1 + e_n u_1 + e_n \underline{\omega} v_2, 
$$
we hence obtain  for all non-vanishing frequencies $\underline{\omega} \neq 0$:
\begin{eqnarray*}
\alpha(\underline{\omega},x_n) &=& x_n^{\frac{k+1}{2}} \Bigg[K_{\frac{k-1}{2}}(\|\underline{\omega}\|x_n)i e_n\Big(-\frac{\underline{\omega}}{\|\underline{\omega}\|} \alpha_2(\underline{\omega})+ \alpha_1(\underline{\omega})\Big)\\ & & \quad\quad\quad\quad  + 
K_{\frac{k+1}{2}}(\|\underline{\omega}\|x_n) \Big(\alpha_2(\underline{\omega})+\frac{\underline{\omega}}{\|\underline{\omega}\|} \alpha_1(\omega)\Big)\Bigg]\\
&=& x_n^{\frac{k+1}{2}}\Bigg[ K_{\frac{k-1}{2}}(\underline{\omega}\|x_n) i e_n \frac{-\underline{\omega}}{\|\underline{\omega}\|}(\alpha_2(\underline{\omega}) + \frac{\underline{\omega}}{\|\underline{\omega}\|} \alpha_1(\underline{\omega}))\\ 
& & \quad\quad\quad\quad  + 
K_{\frac{k+1}{2}}(\|\underline{\omega}\|x_n) (\alpha_2(\underline{\omega}) + \frac{\underline{\omega}}{\|\underline{\omega}\|} \alpha_1(\underline{\omega}))\Bigg]\\
&=& x_n^{\frac{k+1}{2}} \Bigg[K_{\frac{k+1}{2}}(\|\underline{\omega}\| x_n) - i e_n \frac{\underline{\omega}}{\| \underline{\omega}\|}  K_{\frac{k-1}{2}}(\|\underline{\omega}\| x_n)\Bigg] \alpha(\underline{\omega}),
\end{eqnarray*}
where $\alpha(\omega)$ is a Clifford constant belonging to $Cl_{n-1}$. 
\par\medskip\par
It remains to treat the particular case $\underline{\omega} = 0$. In this case, equation (\ref{f1}) simplifies to 
$$
e_n \frac{\partial \alpha}{\partial x_n} + \frac{k}{x_n}(Q \alpha)^{'} = 0.
$$ 
In this case the system (\ref{f1'}) and (\ref{f2'}) simplifies to 
\begin{eqnarray}\label{f7}
\frac{\partial }{\partial x_n} (Q \alpha)^{'} &=& \frac{k}{x_n} (Q \alpha)^{'}\\
\frac{\partial }{\partial x_n} (P \alpha) &=& 0. \label{f8}
\end{eqnarray}
We obtain for $\alpha(\underline{0},x_n) = a(\underline{0}) + \alpha(\underline{0}) x_n^k$, where $a(\underline{0})$ is a constant belonging to the part of the Clifford algebra that has no non-trivial intersection with the vector space $\mathbb{R} e_n$.
\end{proof}
The use of the vector formalism made this proof technically very elegant. Now we can directly translate the resulting  representation formula (\ref{reducedkf}) into the paravector formalism, i.e.,  to the context of  $H^{+}(\mathbb{R} \oplus \mathbb{R}^{n})$:
\begin{corollary}
Suppose that $f:H^+(\mathbb{R} \oplus \mathbb{R}^n) \rightarrow Cl_n$ is a solution to (\ref{kh}). 
Then its image under the Fourier transform in the directions $x_0,x_1,\ldots,x_{n-1}$ have the form 
$$
\alpha(\underline{\omega},x_n) = x_n^{\frac{k+1}{2}} \Bigg[K_{\frac{k+1}{2}}(\|\underline{\omega}\| x_n) - i e_n  \frac{\underline{\omega}}{\| \underline{\omega}\|} K_{\frac{k-1}{2}}(\|\underline{\omega}\| x_n)\Bigg] \alpha(\underline{\omega})
$$
if $\underline{\omega} \in \mathbb{R} \oplus \mathbb{R}^{n-1}\backslash \{0\}$  
and, for $\underline{\omega} \neq \underline{0}$, 
$$
\alpha(\underline{\omega},x_n) a(\underline{0}) + \alpha(\underline{0}) x_n^k
$$
where $a(\underline{0})$ and all $\alpha(\underline{\omega})$ are well-defined Clifford numbers from the Clifford sub algebra $Cl_{n-1}$. 
\end{corollary}
Now suppose that $f$ is a $k$-hypermonogenic automorphic form on $\Gamma_{n-1}({\cal{I}})[N]$. These are in particular  $n$-fold periodic with respect to the discrete period lattice $\Lambda_{n-1}({\cal{I}})[N]$ and we finally arrive at the main result of this subsection: 
\begin{theorem}\label{kf}
Let $f: H^{+}(\mathbb{R} \oplus \mathbb{R}^n) \rightarrow Cl_n$ be a $k$-hypermonogenic automorphic form on $\Gamma_{n-1}({\cal{I}})[N]$.  Then $f$ has a particular Fourier series representation on upper half-space of the form 
{\small{
\begin{eqnarray}
\label{khypfourier}
f(x) &=&  a(\underline{0}) + \alpha(\underline{0}) x_n^k \\
&+&\!\!\!\!\!  \sum_{\underline{m} \in \Lambda_{n-1}^{*}({\cal I})[N] \backslash \{\underline{0}\}} \!\!
x_n^{\frac{k+1}{2}} \Bigg[K_{\frac{k+1}{2}}(2 \pi \|\underline{m}\| x_n) - i e_n  \frac{\underline{m}}{\| \underline{m}\|} K_{\frac{k-1}{2}}(2 \pi \|\underline{m}\| x_n)\Bigg] \alpha(\underline{m}) \nonumber \\
& & \quad\quad\quad\quad\quad \times 
 e^{2 \pi i <\underline{m},\underline{x}>}\nonumber
\end{eqnarray}}}
where $a(\underline{0})$ and all $\alpha(\underline{\omega})$ are well-defined Clifford numbers from the Clifford sub algebra $Cl_{n-1}$. 
\end{theorem}
{\bf Remark.} In the monogenic case $(k=0)$, we have $$K_{-\frac{1}{2}}(\|\underline{\omega}\| x_n) = K_{\frac{1}{2}}(\|\underline{\omega}\| x_n) = \frac{1}{\sqrt{x_n}} e^{-\|\underline{\omega}\|x_n}.$$ 
In this particular case, the factor  simplifies to 
$$
x_n^{\frac{1}{2}}\Bigg[K_{\frac{1}{2}}(\|\underline{\omega}\| x_n) - i e_n \frac{\underline{\omega}}{\| \underline{\omega}\|}  K_{-\frac{1}{2}}(\|\underline{\omega}\| x_n)\Bigg]
= e^{- \pi \|\underline{\omega}\|x_n} (1-i e_n \frac{\underline{\omega}}{\| \underline{\omega}\|}).
$$
This is a scalar multiple of an idempotent in the Clifford algebra. In the monogenic case one obtains the well-known particular Fourier series representation, involving the
monogenic plane wave exponential functions (see \cite{DSS,KraPhD}):
\begin{equation}
\label{fouriermonogenic}
f(x) = \alpha(\underline{0}) + \sum\limits_{\underline{m} \in \Lambda_{n-1}^{*}({\cal{I}})[N] \backslash \{\underline{0}\}}  (1 - i e_n \frac{\underline{m}}{|\underline{m}|}) \alpha(\underline{m}) e^{2 \pi i <\underline{m},\underline{x}> - 2 \pi \|\underline{m}\| x_n}.
\end{equation}
In particular, one has the identity:
$$
\lim\limits_{x_n\rightarrow +\infty} f(x_n e_n) = \alpha(\underline{0}). 
$$ 
\subsection{The Fourier expansion of the $k$-hypermonogenic Eisenstein series}
In this subsection we now apply the previous result in order to determine the Fourier coefficients of the $k$-hypermonogenic Eisenstein series explicitly. To proceed in this direction we first show 
\begin{proposition}
Let $\underline{m} \in \mathbb{R} \oplus \mathbb{R}^{n-1}$ and $x_n > 0$. Then 
\begin{eqnarray} \label{betam}
\beta(\underline{m},x_n) &:=& \int\limits_{\mathbb{R} \oplus \mathbb{R}^{n-1}} \frac{\overline{x}}{\|x\|^{n+1-k}} e^{- 2 \pi i <\underline{m},\underline{x}>} dx_0 dx_1 \ldots dx_{n-1}\\
 & = & \frac{2^{\frac{k+3}{2}} (2 \pi \|\underline{m}\|)^{-\frac{k+1}{2}} \pi^{\frac{n}{2}} x_n^{\frac{k+1}{2}}}{(1-n+k)\Gamma(\frac{n+1-k}{2}-1)} \nonumber \\
 & &  \times \Big(K_{\frac{k+1}{2}}(2 \pi \|\underline{m}\| x_n) - i e_n \frac{\underline{m}}{\|\underline{m}\|}   K_{\frac{k-1}{2}}(2 \pi \|\underline{m}\| x_n)\Big) \nonumber
\end{eqnarray} 
\end{proposition} 
\begin{proof}
We have 
$$
\frac{\overline{x}}{\|x\|^{n+1-k}} = \frac{1}{(1-n+k)} \overline{D} \Big[ \frac{1}{\|x\|^{n-k-1}}\Big].
$$
Now we first compute 
\begin{eqnarray*}
& & \int\limits_{\mathbb{R} \oplus \mathbb{R}^{n-1}} \frac{1}{\|\underline{x} + x_n e_n\|^{n-k-1}} e^{-2 \pi i <\underline{x},\underline{m}>} dx_0 dx_1 \ldots dx_{n-1}\\
&=& \int\limits_{\mathbb{R} \oplus\mathbb{R}^{n-1}} \frac{1}{(x_0^2 + x_1^2 + \cdots + x_{n-1}^2+x_n^2)^{\frac{n-k-1}{2}}} e^{-2 \pi i <\underline{x},\underline{m}>} dx_0 dx_1 \ldots dx_{n-1}\\
&=& \int\limits_{S_{n-2}} \int\limits_{r=0}^{+\infty} \int\limits_{\theta=0}^{+\pi} \frac{r^{n-1}}{(r^2+x_n^2)^{\frac{n-k-1}{2}}} e^{-2 \pi i \|\underline{m}\|\cos(\theta)} \sin^{n-2}(\theta) d\theta dr dS_{n-2},\end{eqnarray*}
where we made the substitution to spherical coordinates by introducing the parameters $r:=x_0^2+\cdots +x_{n-1}^2$ and  $\theta$ as the angle between $\underline{m}$ and $\underline{x}$. $S_{n-2}$ denotes the $n-2$-dimensional unit hypersphere. This integral in turn equals   
\begin{eqnarray*}
& & 2 \frac{\pi^{\frac{n-2}{2}}}{\Gamma(\frac{n-2}{2})} \sqrt{2 \pi}\; 2^{\frac{n+1-k}{2}} \Gamma(\frac{n-1}{2}) (2 \pi \|\underline{m}\|)^{-\frac{n-2}{2}}\\
& & \times \int\limits_{r=0}^{+\infty} \frac{r^{n-1}}{(r^2+x_n^2)^{\frac{n-k-1}{2}}}
r^{-\frac{n-2}{2}} J_{\frac{n-2}{2}}(2 \pi \|\underline{m}\|r) dr\\
&=& \frac{2^{\frac{3+k}{2}} (2 \pi \|\underline{m}\|)^{-\frac{k+1}{2}} \pi^{\frac{n}{2}} x_n^{\frac{k+1}{2}}}{\Gamma(\frac{n+1-k}{2}-1)} K_{\frac{k+1}{2}}(2 \pi \|\underline{m}\| x_n).
\end{eqnarray*}
Applying finally the $\overline{D}$-operator yields the above stated result. 
\end{proof}
Next we observe that we may rewrite the $k$-hypermonogenic Eisenstein series in the form 
$$
E_{k,n-1,N}(x) = 1 + \sum\limits_{M: {\cal{T}}_{n-1}({\cal{I}})[N] \backslash \Gamma_{n-1}({\cal{I}})[N],c \neq 0} \frac{\overline{x} + \overline{d} \;\overline{c}^{-1}}{\|\overline{x}+\overline{d}\;\overline{c}^{-1}\|^{n+1-k}} \frac{\overline{c}}{\|\overline{c}\|^{n+1-k}}.
$$
If $M = \left(\begin{array}{cc} a & b \\ c & d \end{array} \right)$ runs through a subset of representatives of ${\cal{T}}_{n-1}({\cal{I}})[N] \backslash \Gamma_{n-1} ({\cal{I}})[N]$, where we suppose that $N \ge N_0({\cal{I}})$,  then each $r \in \mathbb{Q} + \mathbb{Q} e_1 + \cdots + \mathbb{Q} e_{n-1}$ is represented in the form $c^{-1}d$ exactly once. Then 
$\nu(r) := \frac{\overline{c}}{\|\overline{c}\|^{n+1-k}}$ becomes well-defined, so that we may re-express the Eisenstein series $E_{k,n-1,N}$ in the following form 
$$
E_{k,n-1,N}(x) = 1+ \sum\limits_{r \in \mathbb{Q} \oplus \mathbb{Q}^{n-1}\backslash\{0\}} \frac{\overline{x+r}}{\|x+r\|^{n+1-k}} \nu(r).
$$ From this representation we may readily establish 
\begin{theorem}
Let $k < -1$, $p=n-1$ and $N \ge N_0({\cal{I}})$. Then the Eisenstein series defined in (\ref{khypereisen}) have the Fourier expansion 
\begin{equation}
E_{k,n-1,N} (x) = 1 + \sum\limits_{\underline{m} \in \Lambda_{n-1}^{*}({\cal I})[N]\backslash \{\underline{0}\}} e^{2 \pi i <\underline{m},\underline{x}>} \beta(\underline{m},x_n) \alpha(\underline{m})
\end{equation}
where $\beta(\underline{m},x_n)$ is defined as in (\ref{betam}) and where 
$$
\alpha(\underline{m}) = \sum\limits_{r \;{\rm mod}\; N, r \neq 0} \nu(r) e^{2 \pi i <r,\underline{m}>}
$$
where $r \;{\rm mod}\; N$ indicates that $r$ runs through a set of representatives of $(\mathbb{Q} \oplus \mathbb{Q}^{n-1})/N \Lambda_{n-1}$. 
\end{theorem} 
In the cases where ${\cal{I}}$ is a ring in which each non-zero element has a unique left-sided or right-sided prime factor decomposition, the Clifford group valued expressions $\alpha(\underline{m})$ can directly be related to  Clifford group valued weighted variants of the Epstein zeta function.     

The simplest non-trivial examples for ${\cal{I}}$  with unique left-sided and right-sided prime factor decomposition are for instance the ring of Gaussian integers $\mathbb{Z}[e_1]$, the ring $\mathbb{Z} + \mathbb{Z} e_1 + \mathbb{Z} e_2 + \mathbb{Z} e_3$ (when embedding it into the quaternions by identifying $e_3$ with the $e_1 e_2$) or the Hurwitz order mentioned in Section~2.3.  Notice that the left-sided prime factor decomposition differs in general from the right-sided prime factor decomposition in the non-commutative cases.

Let us write in analogy to \cite{Kri88} ${\cal{R}}(c)$ for a set of representatives of the right cosets $d + c  (N {\cal{I}} \cap \mathbb{R} \oplus \mathbb{R}^{n-1})$. For $c \neq 0$, let $\rho(c)$ be the greatest rational divisor of $c$ in $N {\cal{I}}$, i.e., $\rho(c) := \max_{l \in \mathbb{N}}\{\frac{1}{l} c \in N {\cal{I}}\}$. Adapting from \cite{Kri88}, in this case the expressions $\alpha(\underline{m})$ can be expressed explicitly in the form 
$$
\alpha(\underline{m}) = \sum\limits_{r \in \mathbb{Q} + \mathbb{Q} e_1 + \cdots + \mathbb{Q}e_{n-1} \backslash\{0\} \;{\rm mod}\;N} \nu(r) e^{2 \pi i <r,\underline{m}>}
$$
Next we split the Clifford group valued expressions $\alpha(\underline{m})$ into their real-valued components, i.e. 
$\alpha(\underline{m}) = \sum_{A} \alpha_A(\underline{m}) e_A$. This representation is valued for all orders ${\cal{I}}$, even for those which do not admit a unique left-sided or right-sided prime factor decomposition. If ${\cal{I}}$ however has a unique left-sided or right-sided prime factor decomposition, then each non-vanishing real component $\alpha_A(\underline{m})$ can be represented as follows: 
\begin{eqnarray*} 
\alpha_A(\underline{m}) &=& \Bigg(\sum\limits_{d \in N {\cal{I}}\backslash\{0\}} \frac{d_A \;{\rm sgn}(c_A)}{\|d\|^{n-k}} \Bigg)^{-1} \sum\limits_{c \in N {\cal{I}} \backslash\{0\}} \frac{|c_A|}{\|c\|^{n-k}} \gamma(c,\underline{m}),
\end{eqnarray*}
where 
$$
\gamma(c,\underline{m}) = \left\{ \begin{array}{ll} \rho(c) \|c\|^2 &{\rm\;if}\;<\overline{d} \;\overline{c^{-1}}, \underline{m}> \in  \mathbb{Z} \;\;\;\forall d \in {\cal{R}}(c)\\
0 & {\rm otherwise} \end{array}\right. 
$$
The previous expression then involves divisor sums similarly to those described in \cite{Kri88} for the quaternionic case. 

\par\medskip\par

{\bf Remark.}  
Applying Proposition 1 allows us directly to set up the Fourier expansion for the $k$-hypermonogenic Eisenstein series with positive $k > 1$.

\end{document}